\documentclass[10pt]{amsart}

\usepackage{eucal}
\usepackage{amscd}
\usepackage[all,cmtip]{xy}
\usepackage{rotating}
\usepackage{hyperref,mathrsfs}
\usepackage{tikz}
\usepackage{amsmath,amssymb}
\usepackage{multicol}

\newtheorem{theorem}{Theorem }[section]

\newtheorem{lemma}[theorem]{Lemma}

\newtheorem{corollary}[theorem]{Corollary}
\newtheorem{proposition}[theorem]{Proposition}

\theoremstyle{definition}
\newtheorem{definition}[theorem]{Definition}

\newcommand{\vb}{\subsubsection*{Example.}}

\def\PG{\mathbf{PG}}

\newcommand{\Aut}{\mathrm{Aut}}

\newcommand{\wis}[1]{{\text{\em \usefont{OT1}{cmtt}{m}{n} #1}}}

\newcommand{\A}{\mathbb{A}}

\newcommand{\B}{\overline{\mathbf{B}}}

\newcommand{\Fun}{\mathbb{F}_1}
\newcommand{\Spec}{\wis{Spec}}
\newcommand{\Z}{\mathbb{Z}}
\newcommand{\Proj}{\wis{Proj}}

\newcommand{\mF}{\mathcal{F}}

\newcommand{\fP}{\mathbf{P}}

\newcommand{\mS}{\mathcal{S}}

\newcommand{\mX}{\mathcal{X}}

\newcommand{\F}{\mathbb{F}}

\newcommand{\mY}{\mathcal{Y}}

\newcommand{\PGL}{\mathbf{PGL}}
\newcommand{\PGaL}{\mathbf{P\Gamma L}}

\newcommand{\Sch}{\wis{Sch}}
\newcommand{\bL}{\mathbb{L}}
\usetikzlibrary{matrix}
\usetikzlibrary{arrows}
\usepackage{pgf}
\usepackage{lscape}
\usepackage{listings}

\usepackage{enumitem}

\title{The structure of Deitmar schemes, II. Zeta functions and automorphism groups}
\author{Manuel M\'erida-Angulo and Koen Thas}
\address{Ghent University, Department of Mathematics, Krijgslaan 281, S22 and S25, B-9000 Ghent, Belgium}
\email{manmerang@gmail.com; koen.thas@gmail.com}

\begin{document}
\maketitle

\begin{abstract}
We provide a coherent overview of a number of recent results obtained by the authors in the theory of schemes
defined over the field with one element. Essentially, this theory encompasses the study of a functor which maps certain geometries including graphs to Deitmar schemes with additional structure, as such introducing a new zeta function for graphs.  The functor is then used to determine automorphism groups of the Deitmar schemes and base extensions to fields. \\
\end{abstract}

{\bf R\'{e}sum\'{e}}.\quad 
Nous donnons un vue d'ensemble d'un nombre de r\'{e}sultats r\'{e}cents qui ont \'{e}t\'{e} obtenu par les auteurs dans le domaine de la th\'{e}orie des sch\'{e}mas sur le corps \`{a} un \'{e}l\'{e}ment. Principalement, cette th\'{e}orie concerne lÕ\'{e}tude d'un foncteur qui envoie certains g\'{e}om\'{e}tries (y compris les graphs) sur un sch\'{e}ma de Deitmar avec une structure additionnelle. De cette mani\`{e}re on introduit aussi une nouvelle fonction zeta pour les graphs. Le foncteur est apr\`{e}s utilis\'{e} pour d\'{e}terminer les groups dÕautomorphismes de sch\'{e}mas de Deitmar et de ceux obtenus apr\`{e}s une extension de base dans autres corps.\\



\section{Introduction}

In ``The structure of Deitmar schemes, I'' \cite{KT-Japan} by the second author (Proc. Japan Acad. Ser. A Math. Sci. {\bf 90}, 2014), the author has studied a certain class of Deitmar schemes | which are schemes defined over the field with one element $\Fun$ (cf. \S \ref{Deitmarintro}) | which are naturally associated to what the author called {\em loose graphs}. Loose graphs relax the definition of graph in that edges with $0$ or $1$ vertices are allowed. A simple edge without vertices is a loose graph, for example. One of the main motivations of \cite{KT-Japan} was the fact that affine and projective space Deitmar schemes naturally correspond to loose stars and complete graphs respectively, and hence a natural generalization might lead to a combinatorial, graph theoretical way to study Deitmar schemes and their base extensions to fields. A property that was lacking in \cite{KT-Japan} was 
that if $x$ is a vertex of degree $m$ in a loose graph $\Gamma$, then there is a neighborhood $\Omega$ of $x$ in 
$\mS(\Gamma)$ such that $\mS(\Gamma)_{\vert \Omega}$ is an affine space of dimension $m$; here, $\mS$ denotes the aforementioned association. This makes the association less natural to study (the notion of local dimension is not suited). 

Recently, the authors of the present paper have taken a different turn, and re-defined the map $\mS$ (and called it $\mF$) in order to meet the local dimension property. More details can be found in \S \ref{functprop}. The $\Z$-schemes arising from loose graphs through application of the map $\mF(\cdot)\otimes_{\Spec(\Fun)}\Spec(\Z)$ are of ``$\Fun$-type'' following Kurokawa \cite{Kurozeta}, and thus are provided with a Kurokawa zeta function, as in \cite{Kurozeta}. As such, we can define a new zeta function for (loose) graphs. In \cite{MMKT}, this association is studied in much detail, and we summarize the results in a first part of this paper. Emphasis is put on the mere {\em calculation} of the zeta function, or, in our case, of the corresponding class in the Grothendieck ring of $\Fun$-schemes. This happens through a process called ``surgery,'' which is a stepwise procedure such that in each step an edge with $2$ vertices from a prescribed set of edges, is replaced by two edges with only one vertex. The local dimension rises, and at the end of the process one winds up with a tree. Using precise results for trees, one can then recover the original zeta function.

Although the map $\mF$ is mentioned as a functor in \cite{MMKT}, it is not proven in that paper that it {\em is} a functor. This is done in much details in a second work \cite{MMKT1}. Using functoriality, in {\em loc. cit.} we start a study of automorphism groups of the Deitmar schemes coming from $\mF$ (a more general definition for Deitmar scheme has to be taken into account, as we will explain below). Note that there are several candidates for automorphism groups: one could study combinatorial groups (acting on the underlying incidence geometry), topological ones (acting on the Zariski topological space), projective groups (coming from automorphisms of the ambient projective space), or scheme-theoretic automorphism groups.  Again, we find precise results for trees which suggest a general approach.

In this paper, which should be seen as the natural successor of \cite{KT-Japan}, we want to present the main results of the theory mentioned above. 
Proofs will be published elsewhere.

\section{Deitmar (congruence) schemes}
\label{Deitmarintro}

We consider an ``$\Fun$-ring'' $A$ to be a multiplicative commutative monoid with an extra absorbing element 0.
Let $\Spec(A)$ be the set of all {\em prime ideals} of $A$ together with a Zariski topology. We refer to \cite{Deitmarschemes2} for the definition of prime ideals of a monoid. This topolo\-gical space endowed with a structure sheaf of $\Fun$-rings is called an {\em affine Deitmar scheme}. We define a {\em monoidal space} to be a pair $(X, \mathcal{O}_X)$ where $X$ is a topological space and $\mathcal{O}_X$ is a sheaf of $\Fun$-rings defined over $X$. A {\em Deitmar scheme} is then a monoidal space such that for every point $x\in X$ there exists an open subset $U \subseteq X$ such that $(U, \mathcal{O}_X|_{U})$ is isomorphic to an affine Deitmar scheme.

For a more detailed definition of Deitmar schemes and the structure sheaf of $\Fun$-rings, we refer to \cite{Deitmarschemes2}.

\subsection{Affine space} 

In this paper, $\Fun[X_1, \ldots, X_n]$ denotes the 
{\em monoidal ring} in $n$ variables $X_1, \ldots,\\ X_n$; it is the free abelian monoid generated by $X_1,\ldots,X_n$, containing a 
multiplicative identity $1$ and an absorbing element $0 \ne 1$.

Write $A:=\Fun[X_1, \ldots, X_n]$; then the {\em $n$-dimensional affine space over $\Fun$} is defined as the monoidal space $\Spec(A)$ and denoted by $\mathbb{A}^n_{\Fun}$. All its prime ideals are finite unions of ideals of the form $(X_i)$, where $(X_i)=\{X_ia ~|~ a\in A \}$. 

\subsection{Congruence schemes} 

A more general version of Deitmar scheme is a so-called {\em congruen\-ce scheme} (one can find a more detailed definition in \cite{Deitmarcongruence}), which is defined in terms of {\em sesquiads}. A sesquiad is a monoid $A$ endowed with an {\em addition} or $+$-{\em structure}; this $+$-structure allows addition for a {\em certain} set of elements in the monoid $A$. The category of monoids is a full subcategory of the category of sesquiads.

A sesquiad is said to be {\em integral} if $1\neq 0$, and if from $af=bf$ follows that $(a=b ~ \mbox{ or } ~ f=0)$.

A {\em congruence} on a sesquiad $A$ is an equivalence relation $\mathcal{C}\subseteq A\times A$ such that there is a sesquiad structure on $A/\mathcal{C}$ that makes the projection $A \rightarrow A/\mathcal{C}$ a morphism of sesquiads. If $A/\mathcal{C}$ is integral, the congruence $\mathcal{C}$ is called {\em prime}. We denote by $\Spec_{c}(A)$ the set of all prime congruences on the sesquiad $A$ with the topology generated by all sets of the form 
\begin{equation*}
D(a,b)=\{ \mathcal{C} \in \Spec_{c}(A) ~|~ (a,b)\notin \mathcal{C }\}, ~~~~~~~a,b \in A.
\end{equation*}

In a similar way as for monoids, one can now define a {\em structure sheaf of sesquiads} and a {\em sesquiaded space}. An {\em affine congruence scheme} is a sesquiaded space that is of the form $(\Spec_{c}(A), \mathcal{O}_A)$, for $A$ a sesquiad and $\mathcal{O}_A$ its corresponding structure sheaf, and a {\em congruence scheme} is a sesquiaded space $X$ that locally looks like an affine one. 

\subsection{The $\Proj_{c}$-construction} 

Consider the monoid $\Fun[X_0, X_1, \ldots , X_m]$, where $m\in\mathbb{N}$, as a sesquiad together with the {\em trivial addition}. Since any polynomial is trivially homogeneous in this sesquiad, we have a natural grading
\begin{equation*}
\Fun[X_0, \ldots , X_m]= \bigoplus_{i\geq 0}R_i=\coprod_{i\geq 0}R_i,
\end{equation*}

\noindent where $R_i$ consists of the elements of $\Fun[X_0, \ldots, X_m]$ of total degree $i$, for $i\in \mathbb{N}$. The {\em irrelevant congruence} is given by 
\begin{equation*}
\mbox{Irr}_{c} := \langle X_0 \sim 0 , \ldots, X_m \sim 0 \rangle.
\end{equation*}

Now we can proceed with the usual \Proj-construction of projective schemes. We define\\ 
$\Proj_{c}(\Fun[X_0,\ldots, X_m])$ as the set of prime congruences of the sesquiad $\Fun[X_0,\ldots,X_m]$ which do not contain $\mbox{Irr}_{c}$. The closed sets of the topology are generated by
\begin{equation*}
V(a,b):=\{\mathcal{C}~|~ \mathcal{C}\in \mbox{Proj}_{c}(\Fun[X_0,\ldots, X_m]) , ~~ a\sim_{\mathcal{C}} b\},
\end{equation*}
\noindent for any $(a,b)$ pair of elements of $\Fun[X_0,\ldots,X_m]$. Defining the structure sheaf similarly as in \cite{Deitmarcongruence}, one obtains that $\Proj_c(\Fun[X_0,\ldots,X_m])$ is a {\em projective} congruence scheme. Its closed points naturally correspond to the $\mathbb{F}_2$-rational points of the projective space $\mathbb{P}^{m}(\mathbb{F}_2)$, i.e., elements of 
\begin{equation}
\mathrm{hom}(\Spec(\F_2),\mathbb{P}^m(\F_2)). 
\end{equation}
The space $\mathbb{P}^m(\F_2)$ has a finer subspace structure though, and also a different algebraic structure.

\section{Loose graphs}
A {\em loose graph}\index{loose graph} is a point-line geometry in which each line has at most two different points. Through the analogy with graphs, we call points ``vertices'' and lines ``edges.''
Usually we will consider connected loose graphs, we do not allow loops, and the geometry is undirected. 

Note that any graph is a loose graph.

\subsection{Embedding theorem} 
\label{EmbThm}
 
 Let $\Gamma$ be a loose graph. The embedding theorem of \cite{KT-Japan} observes that $\Gamma$ can be seen as a 
 subgeometry of the combinatorial projective $\Fun$-space $\mathbf{P}(\Gamma)$, called the {\em ambient} space, by simply 
 adding the vertices on each edge which does not contain two vertices, so as to obtain its graph completion, and then constructing the complete graph on the 
 total set of vertices. We will use the same notation $\mathbf{P}(\Gamma)$ for the associated projective space scheme.

\section{Functoriality property}
\label{functprop}

In this section we will briefly describe how one can associate a Deitmar scheme to a {\em loose graph} $\Gamma$ through a functor, which we call $\mathcal{F}$. This functor must obey a set of rules, namely: \medskip

\begin{itemize}
\item[COV]
If $\Gamma \subset \widetilde{\Gamma}$ is a strict inclusion of loose graphs, $\mF(\Gamma)$ also is a proper subscheme
of $\mF(\widetilde{\Gamma})$.
\item[LOC-DIM]
If $x$ is a vertex of degree $m \in \mathbb{N}^\times$ in $\Gamma$, then there is a neighborhood $\Omega$ of $x$ in 
$\mF(\Gamma)$ such that $\mF(\Gamma)_{\vert \Omega}$ is an affine space of dimension $m$.
\item[CO]
If $K_m$ is a sub complete graph on $m$ vertices in $\Gamma$, then $\mF(K_m)$ is a closed sub projective space 
of dimension $m - 1$ in $\mF(\Gamma)$. 
\item[MG]
An edge without vertices should correspond to a multiplicative group.
\end{itemize}
\medskip

Rule (MG) implies that we have to work with a more general version of Deitmar schemes since the multiplicative group $\mathbb{G}_m$ over $\Fun$ is defined  to be isomorphic to
\begin{equation*} \Spec(\Fun[X, Y]/(XY=1)),\end{equation*} where the last equation generates a congruence on the free abelian monoid $\Fun[X,Y]$. (The equation $XY = 1$ is not defined in Deitmar scheme theory.) 
The reader can find a more detailed explanation of this association in \cite{MMKT}.

\begin{theorem}
The map $\mathcal{F}$ is indeed a functor from the category of loose graphs to the cate\-gory of Deitmar congruence schemes. Moreover, for any finite field $k$ (or $\mathbb{Z}$), the lifting map $\mathcal{F}_k(\cdot)=\mathcal{F}(\cdot)\otimes k$ is also a functor.  
\end{theorem}

Let $\Gamma$ be a loose graph and $\mathcal{F}(\Gamma)$ be the Deitmar scheme associated to it. By definition of the functor $\mathcal{F}$, every vertex $v$ of $\Gamma$ defines an affine space over $\Fun$ defined from the ``loose star'' co\-rresponding to $v$. Let us call $v_1, \ldots, v_k$ the vertices of $\Gamma$ and $\Spec(E_i)$ the affine schemes associated to $v_i$, $1\leq i\leq k$. 

\begin{lemma}\label{l1}
For all $1\leq r,s, \leq k,$ $\Spec(E_r)\cap\Spec(E_s)\neq \emptyset$ if and only if $v_r$ and $v_s$ are adjacent vertices in $\Gamma$.
\end{lemma}

\begin{corollary}\label{conF}
Let $\Gamma$ be a loose graph and $\mathcal{F}(\Gamma)$ its loose scheme. Then, $\mathcal{F}(\Gamma)$ is connected if and only if $\Gamma$ is connected.
\end{corollary}

\section{Grothendieck ring of schemes of finite type over $\Fun$}

The $\Spec$-construction on sesquiads (or particularly on monoids with trivial addtion) allow\-s us to have a scheme theory over $\Fun$ defined in an analogous way to the classical scheme theo\-ry over $\Z$. This also allows us to define the {\em Grothendieck ring of schemes over $\Fun$}.

\begin{definition}
The {\em Grothendieck ring of schemes of finite type over $\Fun$}, denoted as $K_0(\Sch_{\Fun})$, is generated by the isomorphism classes of schemes ${X}$ of finite type over $\Fun$, $[X]_{_{\Fun}}$, with the relation
\begin{equation}
[X]_{_{\Fun}}= [X\setminus Y]_{_{\Fun}} + [Y]_{_{\Fun}} 
\end{equation}
for any closed subscheme $Y$ of $X$ and with the product structure given by
\begin{equation}
[X]_{_{\Fun}}\cdot[Y]_{_{\Fun}}= [X\times_{\Fun}Y]_{_{\Fun}}.
\end{equation}
\end{definition}
 
We will later on use the notation $K_0(\Sch_{k})$ for the Grothendieck ring of schemes of finite type over the field $k$, and we will also use the obvious notation $[~\cdot~]_k$. 
We denote by $\underline{\bL}=[\mathbb{A}^1_{\Fun}]_{_{\Fun}}$ the class of the affine line over $\Fun$. Notice that the multiplicative group $\mathbb{G}_m$ satisfies  $[\mathbb{G}_m]_{_{\Fun}}= \underline{\bL} - 1$, since it can be identified with the affine line minus one point.\medskip

\section{Counting polynomial}

Let $\Gamma$ be a loose tree and $\mathcal{F}(\Gamma)$ its corresponding Deitmar (congruence) scheme.  The next result gives us information about the class of $\mathcal{F}(\Gamma)$ in the Grothendieck ring of Deitmar schemes of finite type, $K_0(\Sch_{\Fun})$. We will use the notation $[\Gamma]_{_{\Fun}}$ for the class of $\mathcal{F}(\Gamma)$ in $K_0(\Sch_{\Fun})$ (also when $\Gamma$ is a general loose graph). We adapt the same notation over fields $k$.
\medskip

\begin{theorem}
\label{D3.1}
Let $\Gamma$ be a loose tree. Let $D$ be the set of degrees $\{d_1, \ldots, d_k \}$ of $V(\Gamma)$ such that $1 < d_1 < d_2 < \ldots < d_k$ and let $n_i$ be the number of vertices of $\Gamma$ with degree $d_i$, $1\leq i \leq k$. We call $E$ the number of vertices of $\Gamma$ with degree 1 and $I= \sum_{i=1}^k n_i - 1$. Then the function $\big[~\cdot~\big]_{_{\Fun}}$ is determined as follows:
\begin{equation}
\big[\Gamma\big]_{_{\Fun}} =  \displaystyle\sum_{i = 1}^k n_i\underline{\bL}^{d_i} - I\cdot\underline{\bL} + I + E.
\end{equation}
\end{theorem}


\section{Surgery} 

In order to inductively calculate the Grothen\-dieck polynomial of a $\mathbb{Z}$-scheme
coming from a general loose graph, we introduce a procedure called {\em surgery}. In each step of the procedure we will ``resolve'' an edge, so as to eventually end up with a tree in much higher dimension. One will have to keep track of how the Grothendieck polynomial scheme changes in each step.

\subsection{Resolution of edges}

Let $\Omega = (V,E)$ be a loose graph, and let $e \in E$ have two distinct vertices $v_1, v_2$. The {\em resolution} of $\Omega$
{\em along} $e$, denoted $\Omega_e$, is the loose graph which is obtained from $\Omega$ by deleting $e$, and adding two new
loose edges (each with one vertex) $e_1$ and $e_2$, where $v_i \in e_i$, $i = 1,2$.

One observes that 
\begin{equation}
\mathrm{dim}(\mathbf{P}(\Omega_e)) = \mathrm{dim}(\mathbf{P}(\Omega)) + 2.
\end{equation}

The following theorem reduces the computation of the alteration of the number of $k$-rational points after resolving an edge, to a local problem.

\begin{theorem}[Affection Principle]
\label{AP}
Let $\Gamma$ be a finite connected loose graph, let $xy$ be an edge on the vertices $x$ and $y$, and let $S$ be a subset of the vertex set. 
Let $k$ be any finite field, and consider the $k$-scheme $\mF(\Gamma) \otimes_{\Fun}k$. Then $\cap_{s \in S}\A_s$,where $\A_s$ is the local affine space corresponding to the vertex $s\in S$, changes when one resolves 
the edge $xy$  only if $\cap_{s \in S}\A_s$ is contained in $\fP_{x,y}$, the projective subspace of $\fP(\Gamma) \otimes_{\Fun}k$ generated by $\B(x,1) \cup \B(y,1)$, where $\B(x,1)=\{v\in V(\Gamma) ~|~ d(v,x)\leq1\}$.
\end{theorem}

In terms of Grothendieck classes, we have the following theorem.

\begin{corollary}[Polynomial Affection Principle]
\label{PAP}
Let $\Gamma$ be a finite connected loose graph, let $xy$ be an edge on the vertices $x$ and $y$, let $\Gamma_{xy}$ be the loose graph after resolving the edge $xy$ and let $k$ be any finite field.
Then in $K_0(\texttt{Sch}_k)$ we have
\begin{equation}\label{pap}
[\Gamma]_k - [\Gamma_{xy}]_k = [\Gamma_{\vert \fP_{x,y}}]_k - [{\Gamma_{xy}}_{\vert \fP_{x,y}}]_k.
\end{equation}
\end{corollary}

\medskip
\subsection{Counting polynomial for general loose graphs}

To compute the counting polynomial of a scheme coming from a loose graph $\Gamma$, we proceed as follows: we choose a spanning loose tree $T$ of $\Gamma$ and resolve in $\Gamma$ all edges not belonging to $T$. This yields a loose tree $\overline{T}$ in which we apply the map defined in Theorem \ref{D3.1} to obtain a counting polynomial for $\overline{T}$. Take an edge $e$ now that was resolved and consider the loose graph $\overline{T}_{e}$ in which all other edges except $e$ are resolved, i.e., $\overline{T}_{e}$ is the next-to-last step in the procedure of obtaining $\overline{T}$. Thanks to Corollary \ref{PAP}, we can compute the counting polynomial for $\overline{T}_{e}$ by restricting it to the changes that occur in $\fP_{e}$ (for the concrete formulas of the Affection Principle we refer to \cite[section 11]{MMKT}). By repeating this process as many times as edges were resolved, we can inductively obtain the Grothendieck polynomial of the scheme associated to $\Gamma$. The validity of this process relies on the next theorem.

\begin{proposition} 
Let $\Gamma$ be a loose graph and let $T$ and $\overline{T}$ be defined as above. Then the Grothendieck polynomial in $K_0(\Sch{\Fun})$ of $\mathcal{F}(\overline{T})$ is independent of the choice of the spanning loose tree $T$ of $\Gamma$.
\end{proposition}

\subsection{Lifting $K_0(\Sch_{\Fun})$}
\label{lift1}

In \cite{Deitmarschemes2}, Deitmar explained how one can extend a scheme over $\Fun$ to a scheme over $\Z$ by lifting affine schemes $\Spec(A)$ to $\Spec(A)\otimes_{\Fun}\Z:= \Spec(\Z[A])$, the gluing being defined by the scheme on the $\Fun$-level. The same base extension is also defined for any finite field $k$.  Thanks to the natu\-rality of the base change functor, this lifting is also compatible on the level of the Grothendieck ring of schemes of finite type.

We define $\Omega$ as a linear map from $K_0(\Sch_{\Fun})$ to $K_0(\Sch_k)$, the Grothendieck ring of schemes of finite type over any field $k$,  sending the class $\underline{\bL}$ to $\bL$, the class of the affine line over $k$.

Notice that the function $\Omega$ is then well defined on the subring $\Z[\underline{\bL}]$ of $K_0(\Sch_{\Fun})$.

We denote, from now on, by $[\Gamma]_{_k}$ the class of its lifting $\mathcal{F}(\Gamma)\otimes_{\Fun}k$ in the Grothendieck ring of schemes of finite type over $k$.

\begin{theorem} Let $\Gamma$ be a loose graph. Then $\Omega([\Gamma]_{_{\Fun}})=[\Gamma]_{_{k}}$.\end{theorem}

\section{A new zeta function for (loose) graphs}

Following \cite{Kurozeta}, we say that a $\Z$-scheme $\mY$ is {\em of $\Fun$-type} if its arithmetic zeta function is of the form
\begin{equation}
\zeta_{\mY}(s) = \prod_{k = 0}^m\zeta(s - k)^{\ell_k},
\end{equation}
where $s$ is in $\mathbb{C}$, $m \in \mathbb{N}$, and the $\ell_j$ in $\mathbb{Z}$. (The zeta functions in the right-hand side are Riemann zeta functions.)
Kurokawa then defines the {\em $\Fun$-zeta function} of $\mY$ to be 

\begin{equation}
\zeta_{\mY}^{\Fun}(s) = \prod_{k = 0}^m(s - k)^{-\ell_k}.
\end{equation}

\begin{theorem}
For any loose graph $\Gamma$, the $\Z$-scheme $\chi := \mF(\Gamma)\otimes_{\Fun}\Z$ is of $\Fun$-type. 
\end{theorem}

\begin{definition}[Zeta function for (loose) graphs]
Let $\Gamma$ be a loose graph, and let $\chi := \mF(\Gamma) \otimes_{\Fun}\Z$. Let $P_{\chi}(X) = \sum_{i = 0}^ma_mX^m \in \Z[X]$ be the zeta polynomial obtained after the surgery process (replacing the class $\bL$ by $X)$. We define the {\em $\Fun$-zeta function} of $\Gamma$ as:
\begin{equation}
\zeta^{\Fun}_{\Gamma}(t) := \displaystyle \prod_{k = 0}^m(t - k)^{-a_k}.
\end{equation}
\end{definition}

\vb In the case of a tree, using the notation from Theorem $\ref{D3.1}$, the zeta function is given by
\begin{equation}
\zeta^{\Fun}_{\Gamma}(t)\ \ =\ \ \frac{(t - 1)^I}{t^{E + I}}\cdot\displaystyle \prod_{k = 1}^m(t - k)^{-n_k}.
\end{equation}

\section{Automorphism groups of $\mathcal{F}(\Gamma)$}

Let $\Gamma$ be a loose graph, $\mathcal{F}(\Gamma)$ be its $\Fun$-scheme and $\mathcal{X}_k=\mathcal{F}(\Gamma)\otimes_{\Fun} k$ its extension to a field $k$.

\subsection{Projective automorphism group} 

We define the {\em projective automorphism group} of the scheme $\mathcal{X}_k$, denoted by $\mbox{Aut}^{\mathrm{proj}}(\mathcal{X}_k)$, as the group of automorphisms of the ambient projective space of $\mathcal{X}_k$ stabilizing $\mathcal{X}_k$ setwise, modulo the group of such automorphisms acting trivially on $\mathcal{X}_k$.

\subsection{Combinatorial automorphism group}

We now consider the scheme $\mathcal{X}_k$ as a point-line geometry, where the set of points $\mathcal{P}$ is the set of $k$-rational points of $\mathcal{X}_k$ and the set of lines $\mathcal{L}$ consists of both projective lines (over $k$) and {\em complete affine lines}. A {\em complete affine} line $l$ of $\mathcal{X}_k$ is a line whose projective completion $\bar{l}$ intersects the scheme $\mathcal{X}_k$ in the whole projective line $\bar{l}$ minus one point.
We define the {\em combinatorial automorphism group} of $\mathcal{X}_k$, denoted by $\mbox{Aut}^{\mathrm{comb}}(\mathcal{X}_k)$, as the group of bijective maps $\mathcal{P}\cap\mathcal{L} \rightarrow \mathcal{P}\cap\mathcal{L}$ that preserve incidence.

\subsection{Topological automorphism group}

We define the {\em topological automorphism group} of the scheme $\mathcal{X}_k$, denoted by $\mbox{Aut}^{\mathrm{top}}(\mathcal{X}_k)$, as the group of homeomorphisms of its underlying topological space.
\begin{proposition} 
\label{combtopo}
The combinatorial group of a scheme $\mX_k$ is a subgroup of the topological automorphism group of $\mathcal{X}_k$. 
\end{proposition}

\section{Automorphisms of general loose trees}

Let $T = (V,E)$ be a finite loose tree, and assume its number of vertices is at least $3$. Let $\overline{T}$ be the graph theoretical completion of $T$ | that is, as before, the tree obtained by adding all end points to the edges of $T$. Define the {\em boundary} of $T$, denoted $\partial(T)$, as the set of vertices of degree $1$ in $\overline{T}$. Let $x$ be a vertex which is at distance $1$ from $\partial(T)$ (i.e., is adjacent with at least one vertex of $\partial(T)$). As $\vert V \vert \geq 3$, $x$ is an inner vertex of degree at least $2$.  

Define $k$ and $\mX_k$ as before. Let $\mathbf{PG}(m - 1,k)$ be the ambient projective space of $\mX_k$. Remember that by the embedding theorem, $T$ can be seen as a subgeometry of a projective $\Fun$-space.

Let $I$ be the set of inner vertices of $\overline{T}$, and for any $w \in I$, let $S(w)$ be the subgroup of $\Aut^{\mathrm{proj}}(\mX_k)$ which fixes the $k$-rational points of $\mX_k$ inside all affine subspaces $\widetilde{\A_v}$ (over $k$) which are generated (over $\Fun$) by a vertex $v$ different from $w$ and all directions on $v$ which are not incident with $w$.  For instance, if the distance of $v$ to $w$ is at least $2$, the local space at $v$ is fixed pointwise, and if the distance is $1$, $\widetilde{\A_v}$ is an affine space of dimension one less than the dimension of $\A_v$.
(In particular, the points in $I \cap \mathbf{B}(w,1)$ are fixed.) 

In the next theorem, one recalls that $\mX_k$ comes with an embedding
\begin{equation}
T \ \hookrightarrow \ \mX_k \ \hookrightarrow\ \PG(m - 1,k),
\end{equation}
so that it makes sense to consider stabilizers of substructures of $T$ in, e.g.,  $\PGL(\mX_k)$.

If $S$ is a set of points in $\PG(m-1,k)$, $\PGaL_{m}(k)_{[S]}$ denotes its pointwise stabilizer.\\

\begin{theorem}
\label{cenprod}
Let $\PGL(\mX_k)_{[I]}$ be defined as 
\begin{equation}
\Aut^{\mathrm{proj}}(\mX_k)_{[I]}\ \cap\ \PGL_{m}(k).
\end{equation}
Then $\PGL(\mX_k)_{[I]}$ is isomorphic to the central product
\begin{equation}
\prod^{\mathrm{centr}}_{w \in I}S(w).
\end{equation}
\end{theorem}

\subsubsection{Determination of $S(w)$} 

We distinguish 2 diffe\-rent cases for computing $S(w)$.

\subsection*{\quad $\dagger$ $w$ is the only inner point} 

All the edges are then incident with $w$. 
Call $E$ the set of such edges with an end point, and $L$ the set of loose edges. Put $\vert E \vert = e$ and $\vert L \vert = \ell$. Then
\begin{equation}
S(w) \cong {\Big(\PGaL_{e + \ell + 1}(k)_{L}\Big)}_{[E \cup \{w\}]},
\end{equation}
\noindent since by definition $S(w)$ fixes all end points of edges from $E$.

\subsection*{\quad $\ddagger$ $w$ is not the only inner point} 

Consider another inner vertex $v$ and let us call $W\ne wv$ one of the edges incident with $v$. Then, $S(w)$ must be a subgroup of $\PGL_m(k)$ since the completion of the affine line determined by $(v,W)$ is fixed pointwise.

Let $E$ and $L$ be as before and let $I$ be the 
set of edges incident with $w$ and with another inner point. Put $\vert E \vert = e$, $\vert L \vert = \ell$ and $\vert I \vert = i$. An element of $S(w)$ induces an element of $\PGL(\overline{\A_w})$ (the latter meaning the projective linear group of the local 
projective space at $w$). Moreover, if two elements $\delta, \delta'$ have the same action on $\overline{\A_w}$, the composition $\delta'\delta^{-1}$ fixes all points of $\PG(m-1, k)$. So $S(w)$ is a subgroup of ${\Big(\PGL_{e + \ell + i + 1}(k)_L\Big)}_{[E \cup I \cup \{w\}]}$.

Note that $S(w)$ fixes all inner vertices. It is induced by the projective linear group, in the projective completion, which fixes the projective spaces based at each element of $I$ but ``away from $w$'' pointwise, while additionally fixing $w$ itself.

\subsubsection{Inner Tree Theorem}

The following theorem is a crucial ingredient in the proof of our main theorem for trees. 

\begin{theorem}[Inner Tree Theorem]
\label{innertree}
Let $T$ be a loose tree, and let $k$ be any field. Put $\mX_k = \mF(T) \otimes_{\Fun} k$, and consider the embedding
\begin{equation}
\iota: T \ \hookrightarrow\ \mX_k.
\end{equation}
Let $\Aut(\mX_k)$ be any of the automorphism groups which are considered in this note (i.e., combinatorial, induced by projective space or topological). 
Let $I$ be the set of inner vertices of $T$, and let $T(I)$ be the subtree (not loose) of $T$ induced on $I$. 
Then if $\vert I \vert \geq 2$, we have that $\Aut(\mX_k)$ stabilizes $\iota(T(I))$. Moreover, $\Aut(\iota(T(I)))$ is induced by $\Aut(\mX_k)$.
\end{theorem}

\subsubsection{The general group}

Before proceeding, we need another lemma. We use the notation of the previous subsection.

\begin{lemma}[Field automorphisms]
\label{lemfield}
Let $\PG(m - 1,k)$ be the ambient space of $\mX_k$. 
We have that 
\begin{equation}
{\PGaL_{m}(k)}_{\mX_k}\Big/{\PGL_m(k)}_{\mX_k} \cong k^{\times}.
\end{equation}
\end{lemma}

Using Lemma \ref{lemfield}, the next theorem determines the complete projective automorphism group.

\begin{theorem}[Projective automorphism group]
\label{MTtrees}
Let $T$ be a loose tree, and let $k$ be any field. Put $\mX_k = \mF(T) \otimes_{\Fun} k$, and consider the embedding
\begin{equation}
\iota: T \ \hookrightarrow\ \mX_k.
\end{equation}
Let $I$ be the set of inner vertices of $T$, and let $T(I)$ be the subtree of $T$ induced on $I$. 
We have $\PGaL(\mX_k) = \Aut^{\mathrm{proj}}(\mX_k)$ is isomorphic to 
\begin{equation}
\Big(\Big(\prod^{\mathrm{centr}}_{w \in I}S(w)\Big) \rtimes \Aut(T(I))\Big) \rtimes k^{\times}.
\end{equation}
\end{theorem}

The condition $\vert I \vert \geq 2$ is essential, as the follo\-wing discussion shows.

\subsubsection{The combinatorial automorphism group}

By Theorem \ref{MTtrees}, we can now determine the combinatorial group as well. 

\begin{theorem}[Combinatorial automorphism group]
\label{autcombloose}
Let $T$ be a loose tree, and let $k$ be any field. Put $\mX_k = \mF(T) \otimes_{\Fun} k$, let $I$ be the set of inner vertices, and suppose that $\vert I \vert \geq 2$. Let $\iota$ be as in Theorem \ref{MTtrees}.
Then
\begin{equation}
\Aut^{\mathrm{comb}}(\mX_k) \cong \Aut^{\mathrm{proj}}(\mX_k).
\end{equation}
\end{theorem}

We have seen in Proposition \ref{combtopo} that for each $\mX_k$, the combinatorial automorphism group is a subgroup of the topological automorphism group. One observes that any projectively induced automorphism is combinatorial, but the other direction is in general {\em not} true. 
For example, let $\Gamma$ be an edge with two different vertices, so that for all $k$, $\mX_k$ is a projective $k$-line. Then each permutation of the $k$-points yields a combinatorial automorphism, but not all of these come from projective automorphisms for all $k$. 
So, in general, 
\begin{equation}
\begin{cases}
\Aut^{\mathrm{top}}(\mX_k) \geq \Aut^{\mathrm{comb}}(\mX_k)\\
\Aut^{\mathrm{comb}}(\mX_k), \Aut^{\mathrm{top}}(\mX_k) \geq \Aut^{\mathrm{proj}}(\mX_k).
\end{cases}
\end{equation}


\end{document}